\documentclass[12pt]{elsarticle}
\usepackage{graphicx}
\usepackage{lineno,hyperref}
\usepackage{lscape,float}
\usepackage{xcolor}
\usepackage{amsmath}
\modulolinenumbers[5]

\journal{Operations Research for Health Care}

\begin{document}

\begin{frontmatter}

\title{Long term simulation analysis of deceased donor initiated chains in kidney exchange programs}

\author[muy]{Utkarsh Verma}
\ead{utkarshverma@iitb.ac.in}

\author[muy]{Narayan Rangaraj}
\ead{narayan.rangaraj@iitb.ac.in}

\address[muy]{Department of Industrial Engineering and Operations Research, Indian Institute of Technology, Bombay, India, 400076}


\begin{abstract}
Kidney exchange programs have been developed to find compatible kidneys for recipients who have incompatible donors.  On the other hand, patients who do not have a living donor depend upon deceased donor (DD) donations to get a kidney transplant. Currently in India, a deceased donor donates kidneys directly to a deceased donor wait-list.  The idea of initiating an exchange chain starting from a deceased donor kidney is proposed in a few articles (and recently executed in Italy), but no mathematical formulation has been given for this merger.  We have introduced an integer programming formulation that creates deceased donor initiated chains, considering both paired exchange registry and deceased donor allocations simultaneously. There is a possibility of overlap between paired exchange registry and deceased donor wait-list registry data as recipients can register for both registries independently. This has also been addressed in the paper. A long term simulation study is done to analyze the gain of these DD initiated chains over time. It suggests that even with small numbers of deceased donors, these chains can increases the number of potential transplants significantly. Also, waiting time and dropout rates in the merged registry decreases substantially. 
\end{abstract}


\end{frontmatter}

\section{Introduction}

Transplantation is the preferred treatment over dialysis for patients with kidney failure. The limitation of this therapy is the availability of compatible donors. Waiting time for a compatible kidney is very high and vary across different blood groups [1]. If a recipient has a compatible living donor, then the transplant is performed. Otherwise, an incompatible pair can undergo an ABO-incompatible transplant (if blood group incompatible) or they can register for a Kidney Exchange Program (KEP). Another way to get a compatible kidney is via non directed donors like an altruistic donor or a deceased donor. An altruistic donor is one who donates one of his kidneys to a recipient with kidney failure without receiving any incentive from this act. A Deceased Donor (DD) is one who donates his organs to wait-listed patients in the DD wait-list registry, once he satisfies the criteria of brain stem death. Different countries have different laws that govern these kinds of donations [2]. In countries like India, altruistic donations are not permissible under the law.[3,4]

Kidney Exchange Program (KEP) also known as Paired Kidney Donation (PKE) were initially performed as binary swaps, but later it was recognized that more patients could benefit if such swaps are extended into longer chains or cycles. Although simultaneous execution of very long exchange cycles is challenging for logistical reasons, still cycles of up to 3 to 5 pairs can be successfully performed. One of the critical aspects of implementing long cycles is that all the surgeries in a cycle need to be done simultaneously. It will eliminate the chances of donor backing out. Since if any of the donors refuse to participate in the cycle after his linked recipient receives a kidney, then the whole cycle gets terminated. The subsequent recipient in the cycle will lose his paired donor and also the chance of participating in another exchange at a later stage. Although the non-simultaneous paired exchange is permissible under the law in North America, where long-chain exchanges have been performed, in countries like India, paired exchanges have to be simultaneous. These considerations create a limitation on the cycle length, and most of the KEPs programs include this limitation.

Non-Directed Donor (NDD) chains also become part of the Kidney Exchange Program, here an altruistic or a deceased donor starts a chain of exchanges when he donates a kidney to a pair in the KEP registry [5,6]. Now, the donor of that pair will donate his kidney to the next pair in the chain, and this will go on till a donor fail to donate his kidney.  These transplants are performed in segments, and the donor of the last pair becomes a bridge donor for the next segment. The longest non-directed donor chain performed in the US had more than 100 transplants [7].

Another possibility of exchange is known as List exchange or the indirect exchange. It is an exchange mechanism between a Donor-Recipient (DR) pair and a recipient in the DD wait-list. Here wherein the donor of an incompatible pair donates his kidney to a recipient in the DD wait-list and in return, his intended recipient gets a high priority in that DD wait-list. It improves the chances of the recipient of that pair to get a compatible DD kidney sooner as compared to a long waiting time for a usual DD kidney. Indirect exchange works as follows, an incompatible pair donates his kidney to a high-rank DD wait-list recipient and in return recipient of the pair receives the rank of DD recipient. However, this process also has challenges in implementation because the arrival of DD kidney is random. A Paired recipient may have to wait for a long time to get a compatible kidney, and in the meantime, the patient might become too ill for transplant or even die, so even after donating his donor's kidney, he may not get a compatible kidney. Also, if a patient donates his intended donor's kidney, then he will not be able to take part in possible future KEP and will lose the chance of getting a living donor. These things make indirect exchange practically hard to implement.

Currently deceased donor allocations are done independently of kidney exchange program and an idea of initiating deceased donor chains were proposed in the literature [3,8]. A two length deceased donor chain was performed in Italy last year. A study of creating deceased donor chains and some retrospective analysis of possibilities of such chains was done recently [9]. Extending the analysis of such chains, a long term simulation study of deceased donor chains is done here on some realistic data from Indian registries. In this paper, we have analyzed the deceased donor initiated chains and performed long term simulations to observe the effect of these chains on KEPs. 

The major contributions of this paper are as follows:

\begin{enumerate}
    \item An Integer Programming model for creating DD initiated chains is proposed.
    \item Recipients can register for KEP as well as DD wait-lists, and this overlap is explicitly addressed in the proposed model.
    \item Long term simulation analyses of deceased donor initiated chains are performed, and their comparison with the current process is shown.
    \item Performance measures such as total number of transplants, waiting time and dropout rates are also compared.
    \item The numerical simulations are tested for realistic data from live KEP registries and representative blood group data for deceased donors (in India).
\end{enumerate}

In the next section, a brief literature review of kidney exchange programs is given, and later, a mathematical model and analysis of deceased donor initiated chains are shown.

\section{Literature Review}

Kidney exchange programs (KEP) or paired kidney donation (PKD) have emerged over the last two decades [10-16]. It was first introduced by Rapaport in 1986 when he proposed an idea of creating a living donor pool for kidney exchange [17]. In 2000, the US initiated a pilot program for paired kidney exchange and later other countries like the UK, Netherlands and Italy started their kidney exchange program [10,11]. Initially, the pool size of these exchange programs was small, and sometimes the proposed transplants also fails due to tissue type incompatibility [18]. Researchers have tried to increase the supply of donors by other mechanisms. Indirect exchange or List exchange was another mechanism where an incompatible pair will trade his donor kidney for a rank in the deceased donor wait-list [12]. These had some ethical concerns as deceased donor arrivals are random, and by donating his donor's kidney, a pair loses the possibility of further direct exchanges. ABO-incompatible or desensitized transplants can also be part of KEP where incompatible pairs can exchange their donors to improve their compatibility [19]. 

Non-directed donor (NDD) initiated chains also became part of KEP, and several mechanisms to create non-directed donor initiated chains were proposed [20-21]. There are two types of non directed donors, deceased donors and altruistic donors. Altruistic donations give an opportunity to create long exchange chain, because these chains can be performed in segments, and the last donor of each segment becomes a bridge donor which can start another segment, and this chain can continue forever in theory [22].

Deceased donor chains are another possibility of exchange chains discussed in recent years. Melcher et al. proposed non-simultaneous altruistic chains starting from deceased donors [8]. These DD initiated chains would benefit DD wait-listed patients as patients of KEP also registers for the deceased donor's kidney. If KEP patients get a compatible kidney through chains, then it will clear the wait-list faster. It can also increase the welfare of deceased donor patients of the chain since they will get a living donor kidney, which has greater graft survival. However, they proposed it as a non-simultaneous chain which can create issues of fairness. If a chain is not executed fully due to donor backing out, then a deceased donor wait-list patient won’t get a compatible kidney which he would have got if the current process of allocation was performed. 

Billa et al. proposed the idea of simultaneous execution of DD initiated chains where multiple alternatives can be created to increase the probability of successful execution [3]. Non-simultaneous executions are not permitted under the law in India. They also showed that there would not be any loss in terms of the number of recipients matched for different blood groups, if all KEP recipients also register for a DD wait-list. 

The first successful execution of the deceased donor initiated chain was reported in Italy last year [23]. This suggests that deceased donor chains might be one way forward in kidney transplants as this will increase the possibility of larger exchanges. There are some ethical concerns regarding DD chains such as unfairness towards type O recipients in DD wait-list, the quality difference between DD and donor provided to DD wait-list, acceptance of DD kidney in return of a living donor kidney and donor backing out from chain which will reduce a transplant in DD wait-list [8]. Rees et al. showed that these chains would not disadvantage O type and highly sensitized patients [24], but there could be other concerns in the implementation of DD chains. In the next section, we extended the idea of Billa et al. [3] and provided an integer programming formulation and long term simulation analysis of the simultaneous execution of DD chains.

\section{Deceased donor initiated chains}
In the currently existing scheme of organ allocation in India, both kidneys of a deceased donor are allotted to the DD wait-list recipients. If one of the kidneys of the DD is used to start a chain exchange then a few more pairs will get a compatible kidney. We propose that the first deceased donor kidney goes to a recipient in the DD wait list registry to maintain the fairness criteria of the recipients on DD wait-list. The other kidney goes to the KEP registry and the paired donor of that recipient donates his kidney to the next recipient in the same KEP registry, to create a paired exchange chain. Several such donor-recipient pairs can then be added to this chain. The last donor of the pair then donates his kidney to the next patient in the deceased donor wait list. Thus, from a single deceased donor, two patients in the DD wait-list receive kidneys while some more patients in the paired exchange registry also simultaneously receive compatible kidneys. 

\begin{figure}
    \centering
    \includegraphics[scale=0.3]{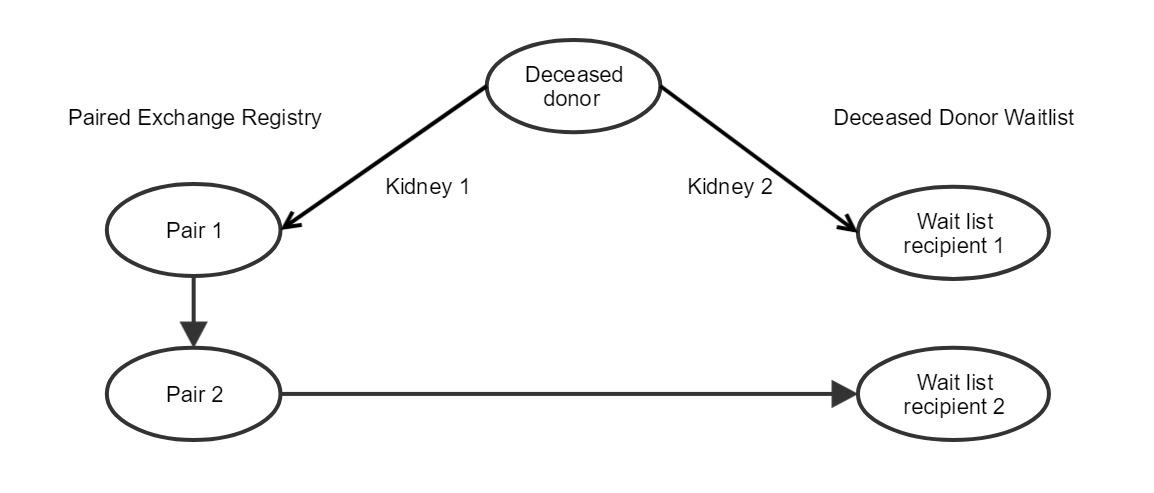}
    \caption{Deceased donor initiated chain of length 3 [3]}
    \label{fig:my_label}
\end{figure}

The aim of creating DD chains is to maximize the utility of a deceased donor kidney. Currently, a deceased donor provides two kidneys to DD wait-list, if we can create DD chains combining KEP registry then a few additional pairs will get a compatible kidney. There will not be any loss to the DD wait-list as the last donor of the chain donates a kidney back to DD wait-list. These chains need to be executed simultaneously to avoid the chances of non-execution of chains, so a bound on chain length is required. We considered 2 length chains, which seem to be feasible after discussing with the medical community. Also, DD initiated chains increase the number of transplants and thus the waiting time for getting a compatible kidney in KEP registry will decrease as the flow of kidney will increase in that pool. In the next section, we will discuss an Integer Programming model for DD initiated chains.

\subsection{IP model for DD initiated chains}

Consider a directed graph $G=(V,E)$ with set of nodes $V$ consisting of four subsets of nodes, set of deceased donors (DD), set of incompatible DR pairs who only register for living donor kidney exchanges(P), set of all incompatible DR pairs who register for both registries (PWL) and set of all DD wait-list recipients (WL). The edge set E consists of all directed edges in the graph and an edge is formed from one node to another if the donor at first node is compatible with the recipient at the other node. So each node in DD will have only outgoing edges from it, each node in WL will have only incoming edges. Nodes in P and PWL will have both incoming and outgoing edges. A weight set $w$ is also formed consisting of weight for each edge where weights are obtained by considering factors like HLA mismatch, age difference, waiting time on dialysis, ranking in wait-list etc. Since solving a bounded cycle and chain KEP is NP hard problem, a compact formulation was created so that the running time of the model is reasonable with a large number of nodes and edges. For compact formulation, L copies of the graph were created where L is the sum of number of DD and half of the number of pair nodes.

Let the variables be

\begin{eqnarray}
x_{i,j}^{l} = 
\left\{
\begin{array}{llll}
1 &  \text{If edge $(i,j)$ is selected in the $l^{th}$ copy of the graph}\\
0 & otherwise
\end{array}
\right.   
\label{DDIC variable}
\end{eqnarray}

$f_{i}^{l}$ - Sum of all outgoing edges for a node $i$ in $l^{th}$ copy of graph 

$g_{i}^{l}$ - Sum of all incoming edges to a node $i$ in $l^{th}$ copy of graph.\\
\\
An integer programming formulation for DD initialed chains is as follows:\\
\\
\begin{equation}
max \sum_{l \in L} \sum_{(i,j)\in E} w_{(i,j)}^{l} \  x_{(i,j)}^{l}
\end{equation}

\begin{equation}
s.t. \hspace{1cm} \sum_{j:(i,j) \in E} x_{(i,j)}^{l} = f_{i}^{l} \hspace{1cm} \forall i \in N, \forall l \in L
\end{equation}

\begin{equation}
 \sum_{j:(j,i) \in E} x_{(j,i)}^{l} = g_{i}^{l} \hspace{1cm} \forall i \in N, \forall l \in L
\end{equation}

\begin{equation}
f_{i}^{l} = g_{i}^{l} \hspace{1cm} \forall i \in P, \forall l \in L
\end{equation}

\begin{equation}
g_{i}^{l} \leq 1 \ \ \forall i \in P,\forall i \in PWL, \forall i \in WL, \forall l \in L
\end{equation}

\begin{equation}
f_{i}^{l} \leq 1 \hspace{1cm} \forall i \in DD,\forall i \in PWL \forall l \in L
\end{equation}

\begin{equation}
f_{i}^{l} \leq g_{i}^{l} \hspace{1cm} \forall i \in PWL, \forall l \in L
\end{equation}

\begin{equation}
\sum_{(i,j) \in E} x_{(i,j)}^{l} \leq k \hspace{1cm} \forall l \in L
\end{equation}

\begin{equation}
\sum_{l \in L} f_{i}^{l} \leq 1 \hspace{1cm} \forall i \in N
\end{equation}

\begin{equation}
\sum_{l \in L} g_{i}^{l} \leq 1 \hspace{1cm} \forall i \in N
\end{equation}

\begin{equation}
x_{ij}^{l} \in \{0,1\}. \hspace{1cm} \forall (i,j) \in E,  \forall l \in (1 ... L)
\end{equation}

Here the objective is to maximize the weighted sum of the total number of edges over all replications. Each edge represents a compatible match and weights are calculated based on a scoring system. Thus we are maximizing the number of feasible transplants while considering the quality of a match.

Constraint (4) ensures that for all pair nodes, sum of incoming edges must be equal to sum of outgoing edges in each replication. It means a pair will only donate his donor's kidney if his intended recipient receives a compatible kidney in that replication else it won't donate. Constraint (5) ensures that maximum 1 kidney can be received in any given replication. Constraint (6) says that for each deceased donor can donate maximum of 1 kidney in any given replication.

Now for those recipients who are on both DD wait-list and swap registries, a new constraints was introduced. Constraint (8) ensures that sum of all outgoing edges should be less than the sum of all incoming edges in a replication for a PWL pair, this will mean that whenever a PWL pair donates his donor's kidney then the recipient also gets a compatible kidney in the same replication but it is not necessary vice-verse. Note here that since these pairs are registered for both registries thus they should get chance to receive a compatible kidney either via a living donor or via a DD kidney depending upon his ranks in both registries. This was done by separating out those pairs who registered for both registries and defined a new variable for them. It ensures that there is no loss to recipients who registers for both registries in the merging registry.

Constraint (8) ensures that maximum donation in a replication is bounded by k. Thus a chain or cycle can be formed with a maximum length of k. Constraint (9) ensures that for a node $i$, maximum donation over all replications can maximum be 1 and (12) ensures that maximum 1 kidney can be received to a recipient over all replications. Constraint (13) represents that variable $x_{(i,j)}^{l}$ is binary variable.

This IP model can be considered as a combination and extension of earlier discussed models by Constantino et. al [25] and Anderson's et. al [26]. The major difference between the proposed model and earlier models is the inclusion of different types of nodes and constraints. Constantino's model is designed for cycle formulation while Anderson's model talks about finding long chains initiated by altruistic donors. We have proposed a model considering both PKE and DD wait-list registry simultaneously. This creates a possibility of overlap between both registries, so we considered this possibility as well and defined a new variable which represents a pair in both registry. Our model has a non-directed donor node (Deceased donor node) which doesn't have any incoming edges and DD wait list node which doesn't have any outgoing edges. A Few additional constraints were also added to create a bound on chain and cycle length as we are considering simultaneous executions of these exchange transplants.

Melcher's et al. proposed that these DD initiated chains can be non-simultaneous and every DD kidney can be considered as a Non directed donor's kidney [8]. But this will create fairness issues as top wait-listed recipients in DD wait-list registry may have to wait long to get a compatible kidney. We have proposed that one of the two DD kidneys should be given directly to the DD wait-list to maintain the fairness criteria and only one kidney should be used to start DD chain [9]. Also, non-simultaneous executions may lead to breaking down of a chain and the DD wait-list won't get a kidney which they would have got otherwise. So to maintain the fairness we propose the executions of the chain to be simultaneous.

\subsection{Different variants of the problem}

Here two variants of the DD initiated chains are discussed which includes patients with multiple incompatible donors and patients with a low matched compatible donor.
\subsubsection{DD initiated chains with Multiple donors}
It can happen that a patient has more than one willing but incompatible donors, these donors can be of different blood groups and ages. So having more than one willing donor will increase the chances of a recipient to get a compatible kidney faster. For these patients, the network will be created considering all willing donors and edge weights will be calculated for each edge. Now when the program runs it consider all edges and tries to maximize the overall weight. The edge which contributes more to the solution will be selected. For example, if the number of recipients benefited is same for two solutions then the edge which has the higher overall weight will be considered and the donor corresponding to that edge will be selected.
\subsubsection{DD initiated chains with compatible pairs}
Compatible pairs can also be part of a DD initiated chain. Pair with low HLA match or large age difference or HIV positivity between donor and recipient, participates in KEP to increase the post-transplant benefits of the recipient. A Few compatible pairs with O donor also join the registry to increase the welfare of kidney recipients since O recipients can receive kidney only from an O donor.\\
For these pairs, an additional constraint can be included which will ensure that the recipient of the compatible pair will get a better-matched kidney than his own match.
\begin{equation}
w_{(j,i)}^{l} \  x_{(j,i)}^{l} \geq w_{(i,i)}^{l} \  x_{(j,i)}^{l}
\end{equation}
\begin{equation}
\forall (j,i) \in E, \forall i \in CN \ \forall l \in (1 ... L) 
\end{equation}

Where CN represents the set of compatible pair nodes. This constraint will ensure that in each replication, the node of the compatible pair should have an incoming weight higher than his own edge weight.

\section{Comparison between different types of kidney exchange programs}

There are several extensions of kidney exchange program developed over the years. A comparison of different extensions of KEP with our model is shown in table 1. This will give us a better picture of where our model currently exists. The table has compared the extensions of KEP, major results and scale of the problem covered.

\begin{table}
    
\caption{Comparison among different types of models proposed for kidney exchange program in literature}
\resizebox{\textwidth}{!}{%
\begin{tabular}{p{5cm}p{5cm}p{5cm}p{3cm}}
\hline
Title \& Reference & Idea proposed & Major results & Scale of the problem\\
\hline
Kidney exchange by Roth et. al. 2004 [12] & Indirect exchange and Top Trade Cycle and Chain(TTCC) algorithm & TTCC can increase the welfare of kidney recipients & Tested upto 300 pairs\\
\hline
Clearing Algorithms for Barter Exchange Markets: Enabling Nationwide Kidney Exchanges by Abraham et. al 2007 [27] & ILP formulation and incremental formulation approach, Cycle of length k is considered & KEP with length more than 2 is NP-complete, Incremental formulation approach performs better than existing approaches & Solved KEP upto 10000 pairs\\
\hline
New insights on integer-programming models for the kidney exchange problem by Constantino et. al. 2013 [25] & New compact formulations, Egde assignment and extended edge formulation & Compact formulations have advantages over non-compact
ones when the problem size is large & Tested upto 1000 pairs with bound on cycle length upto 6\\
\hline
Finding long chains in kidney exchange using the
traveling salesman problem by Anderson et. al. 2014 [26] & Recursive Algorithm for the KEP, PC-TSP–Based Algorithm for the KEP & Long chains can increase the benefit of KEP  & Number of pairs considered upto 350 and number of non deceased donor upto 10\\ 
\hline
IP Solutions for International Kidney Exchange
Programmes by Peter Biro et al. 2019 [28] & IP model for international exchanges with country specific constraints & Merged pool performed better than individual allocations & 20 instances of 1000 pairs were considered\\
\hline
Optimal integration of desensitization protocols into kidney paired
donation (KPD) programs by Fatemeh Karami et al. 2019 [19] & Inclusion of desensetized transplants (ABO incompatible transplants) in KEP & Even with limited number of ABO incompatible transplants, an increase of 66\% can be observed in transplant rates& 826 pairs were considered for simulations \\
\hline
Deceased donor initiated chains in kidney exchange program (proposed model) & IP model for deceased donor initiated chains is proposed considering DD wait list recipients as well & Even with small DD chains, a significant benefit can be achieved through it & For 2 way chains - upto 500 pairs\\
\hline
\end{tabular}
}
\end{table}

\section{Simulation Results}

To do the simulations, we started with the smallest possible DDIC which would be of length 2. We compared the long term behaviour of DDIC over independent functioning of KEP and DD allocations (Current allocation process). The edge weights were considered to be 1 for each edge. The KEP data distribution was taken from ASTRA registry in India and deceased donor blood group distribution was taken from the general blood group distribution in India.

\subsection{Long term analysis of DDIC}

We simulated our model over several months and compared the behaviour of both allocation processes. It requires several assumptions which are as follows:

\begin{enumerate}
  \item Blood Group distribution of pairs follow ASTRA registry data
  \item Blood Group distribution of Deceased donors will follow general Blood Group distribution in India
  \item We will simulate our model for a period of 5 years and KEP runs every month and DDIC are created whenever they become available (60 rounds of simulations)
  \item We replicated our simulation 30 times to average out the randomness of the process.
  \item Three instances of arrival rates were considered for KEP registry,  Uniform(10,15), Uniform(15,20) and Uniform (20,25).
  \item Three instances of DD arrival rates were considered, Uniform(1,5), Uniform(5,10) and Uniform(10,15).
  \item It is assumed that there is a large wait-list for DD kidneys so all the kidneys offered to them will be consumed.
  \item Probability of failure for a un-match pair in KEP registry to go to next round was considered to be, 0(No dropout), 0.1 and 0.3.
\end{enumerate}

The comparison of total number of recipients matched and total number of dropouts was done in each round for both process. In Figure 2, comparisons are shown with different arrival rates and dropout probabilities. It can be observed that, DDIC outperformed the current allocation process both in terms total number of transplants and total number of dropouts. Also, DDIC performed better than the current process in each round as they had higher number of transplants and lesser number of dropouts. 

It was expected that with increase in DD arrival rates, benefit of DDIC will increase proportionally . With higher arrival rates, the graph shows a significant benefit in each round thus it will have a greater impact over all. This happens due to increase in kidney supply as KEP data will have very few O type donors and AB type recipients.

\begin{figure}[H]
\hspace{-2.5cm}
    \begin{minipage}[b]{0.6\textwidth}
		\includegraphics[scale=0.3]{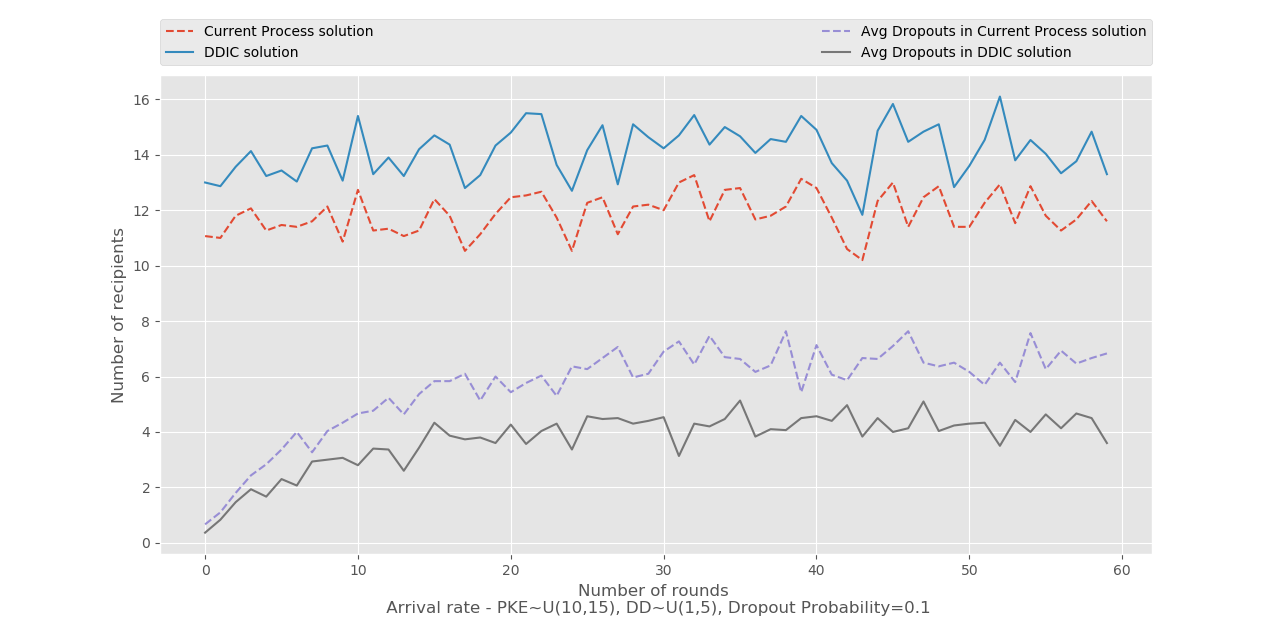}
	\end{minipage}
	\hspace{0.5cm}
	\begin{minipage}[b]{0.6\textwidth}
		\includegraphics[scale=0.3]{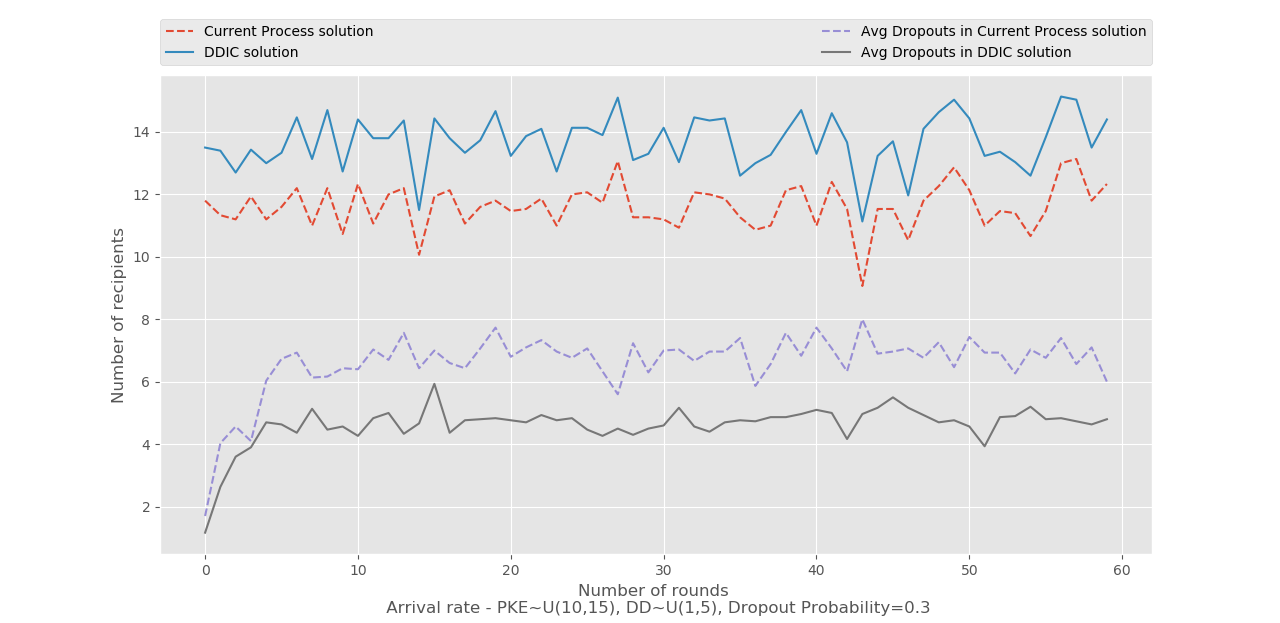}
	\end{minipage}
	
\hspace{-2.5cm}	\begin{minipage}[b]{0.6\textwidth}
		\includegraphics[scale=0.3]{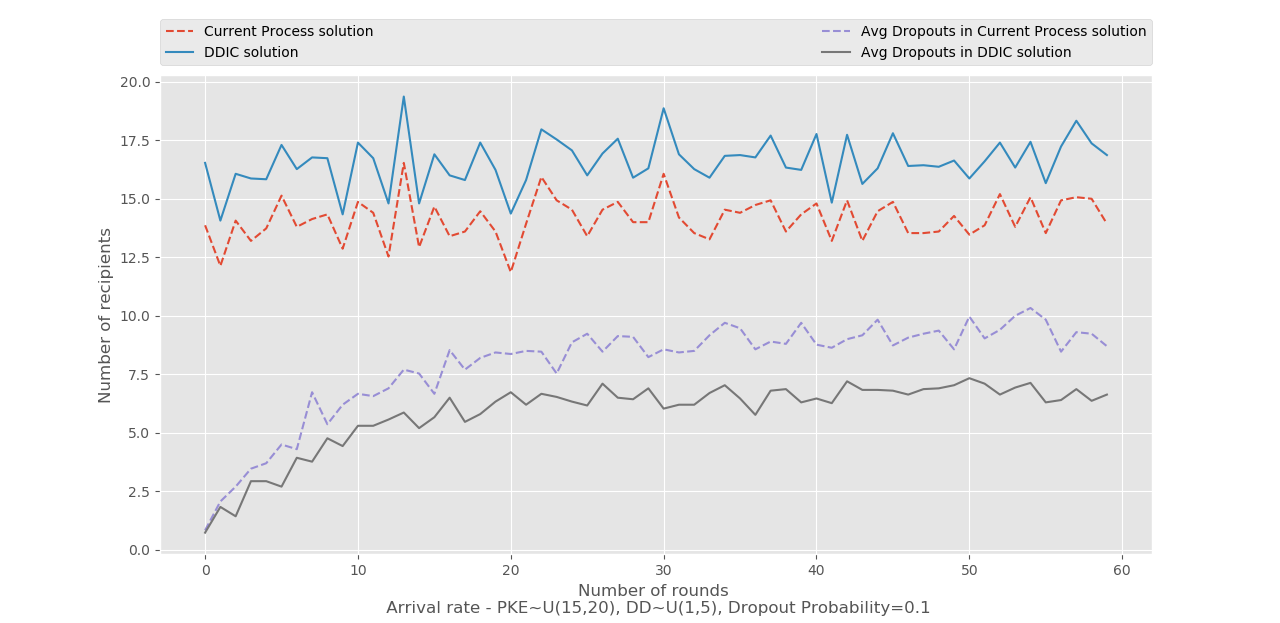}
	\end{minipage}
	\hspace{0.5cm}
	\begin{minipage}[b]{0.6\textwidth}
		\includegraphics[scale=0.3]{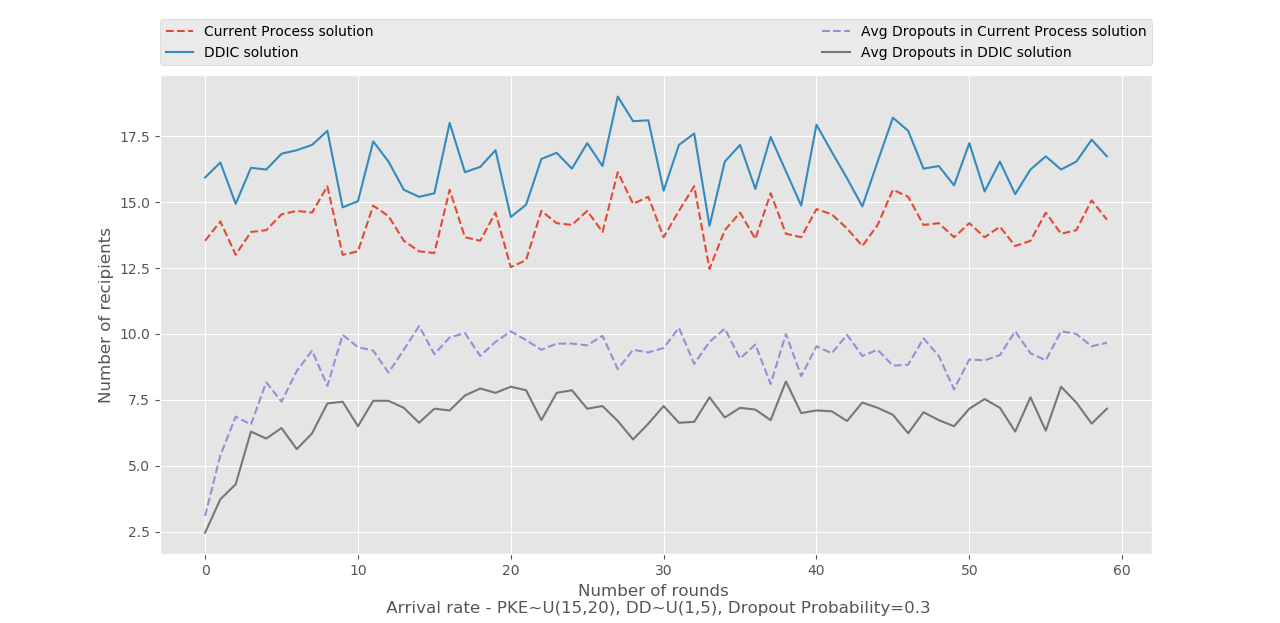}
	\end{minipage}
	
\hspace{-2.5cm}	\begin{minipage}[b]{0.6\textwidth}
		\includegraphics[scale=0.3]{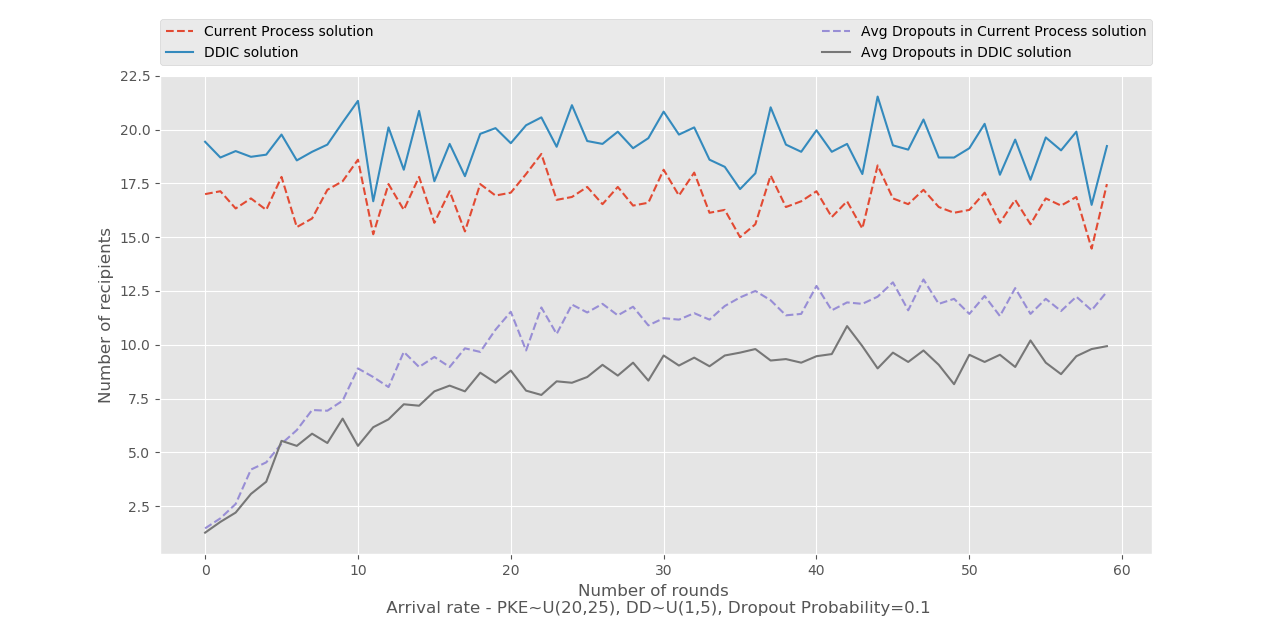}
	\end{minipage}
	\hspace{0.5cm}
	\begin{minipage}[b]{0.6\textwidth}
		\includegraphics[scale=0.3]{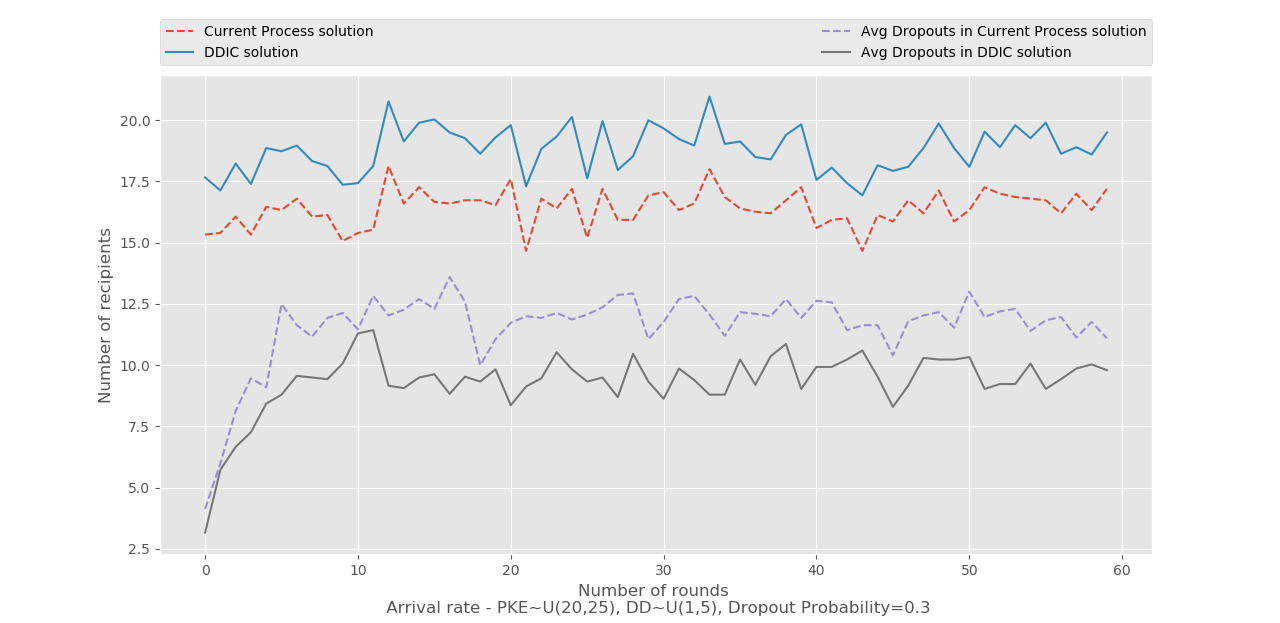}
	\end{minipage}
\end{figure}
\begin{figure}[H]
\hspace{-2.5cm}	\begin{minipage}[b]{0.6\textwidth}
		\includegraphics[scale=0.3]{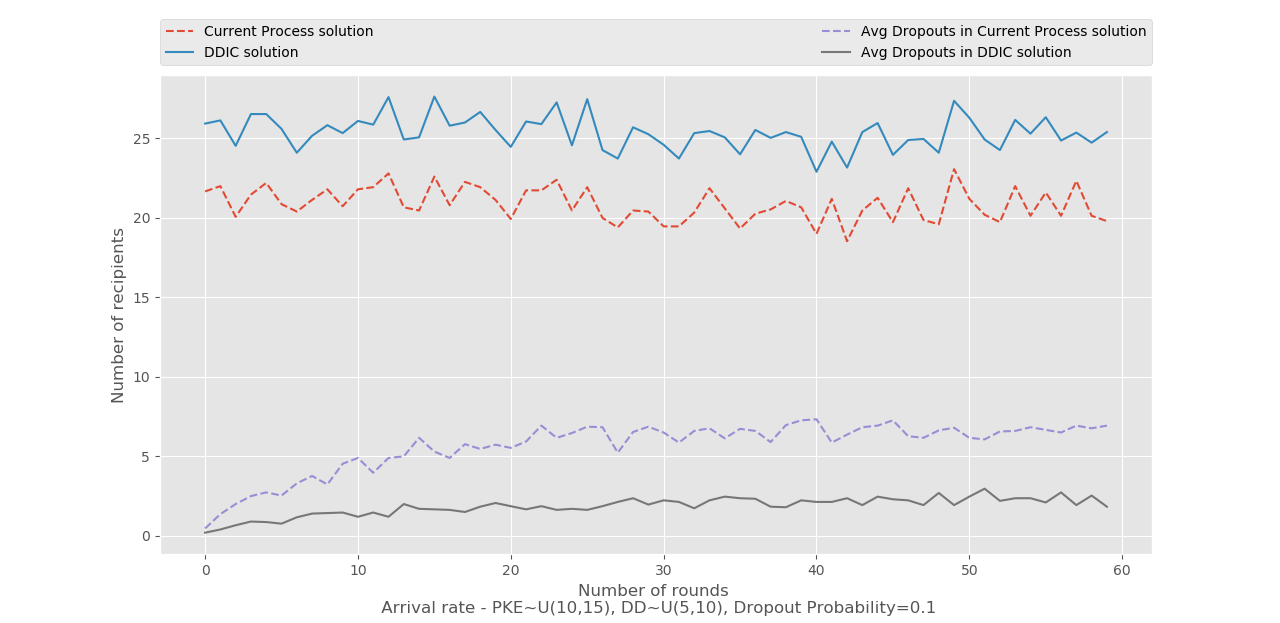}
	\end{minipage}
	\hspace{0.5cm}
	\begin{minipage}[b]{0.6\textwidth}
		\includegraphics[scale=0.3]{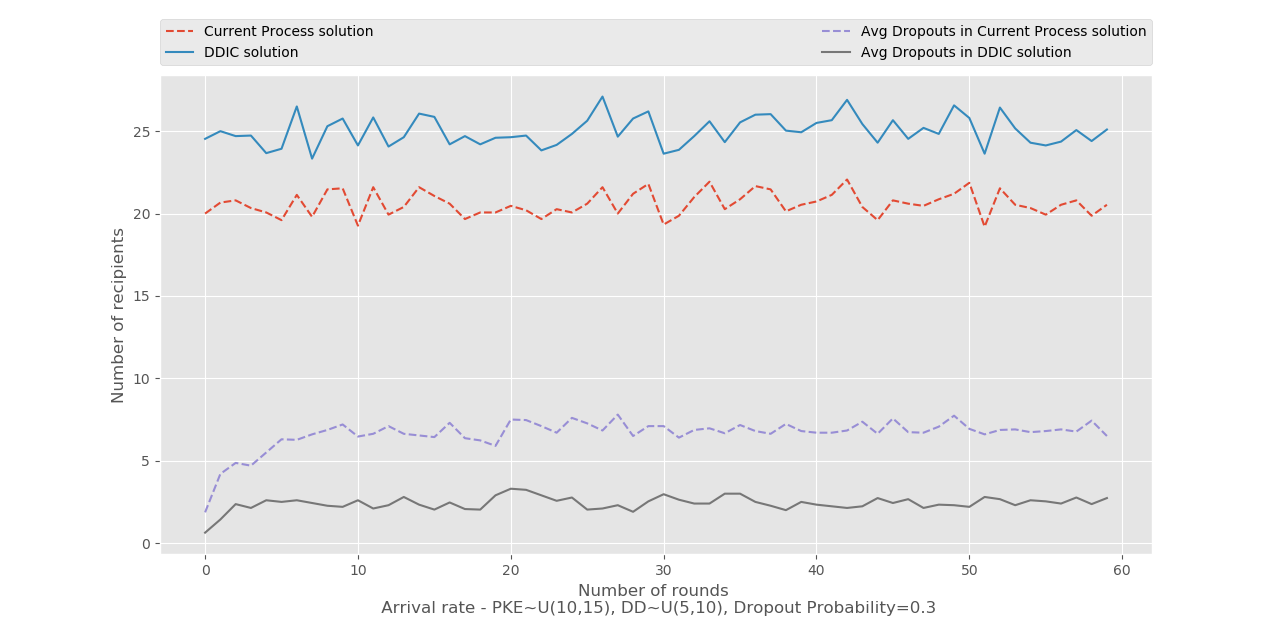}
	\end{minipage}
	
\hspace{-2.5cm}	\begin{minipage}[b]{0.6\textwidth}
		\includegraphics[scale=0.3]{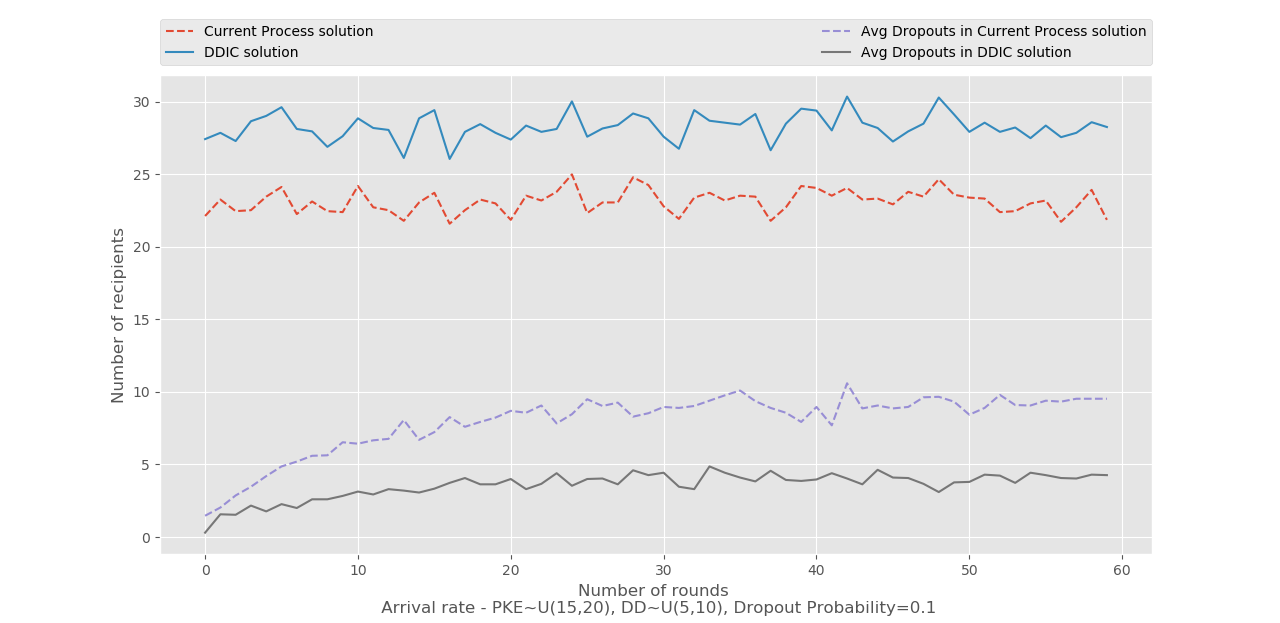}
	\end{minipage}
	\hspace{0.5cm}
	\begin{minipage}[b]{0.6\textwidth}
		\includegraphics[scale=0.3]{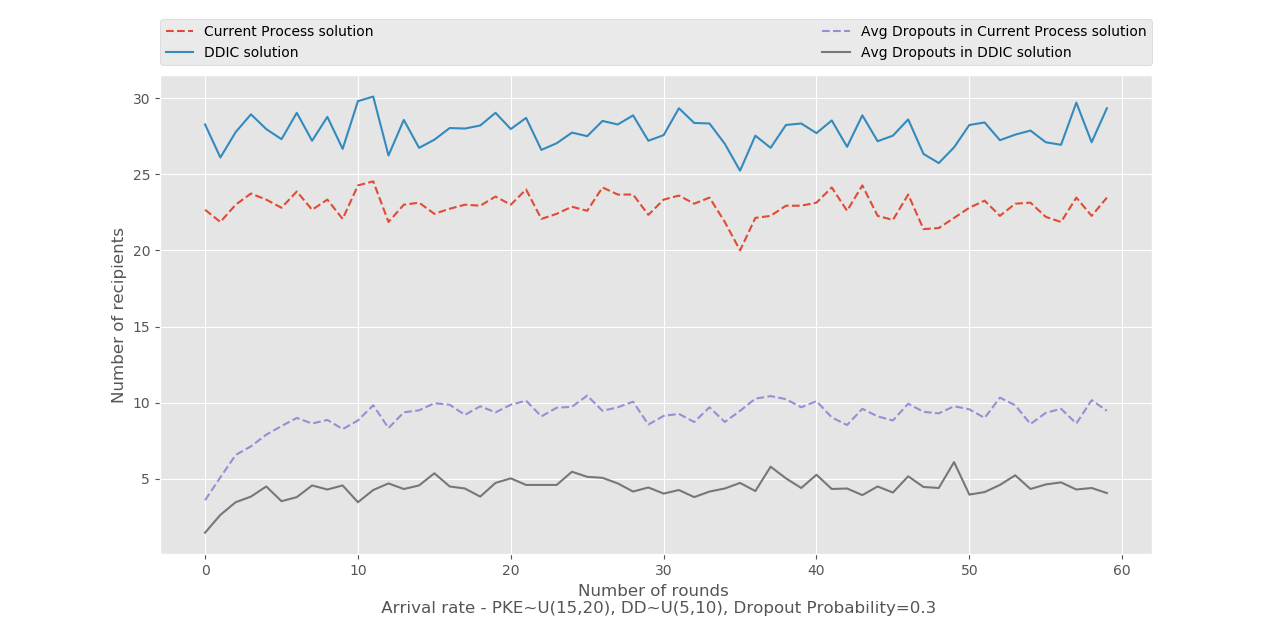}
	\end{minipage}
	
\hspace{-2.5cm}	\begin{minipage}[b]{0.6\textwidth}
		\includegraphics[scale=0.3]{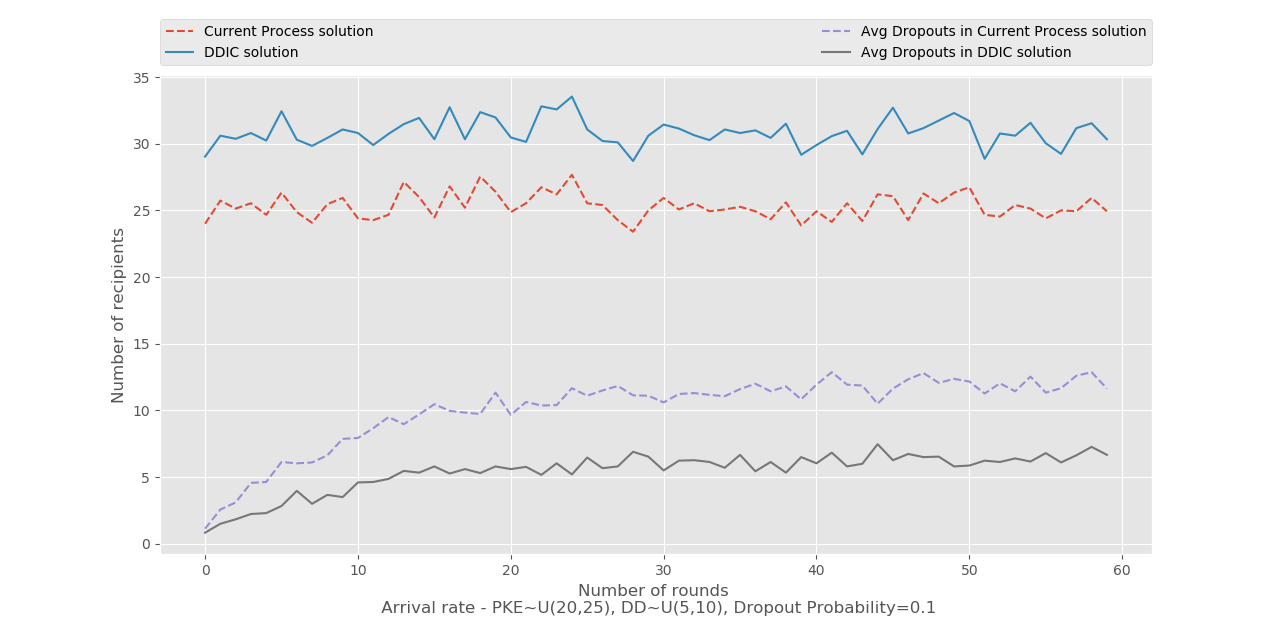}
	\end{minipage}
	\hspace{0.5cm}
	\begin{minipage}[b]{0.6\textwidth}
		\includegraphics[scale=0.3]{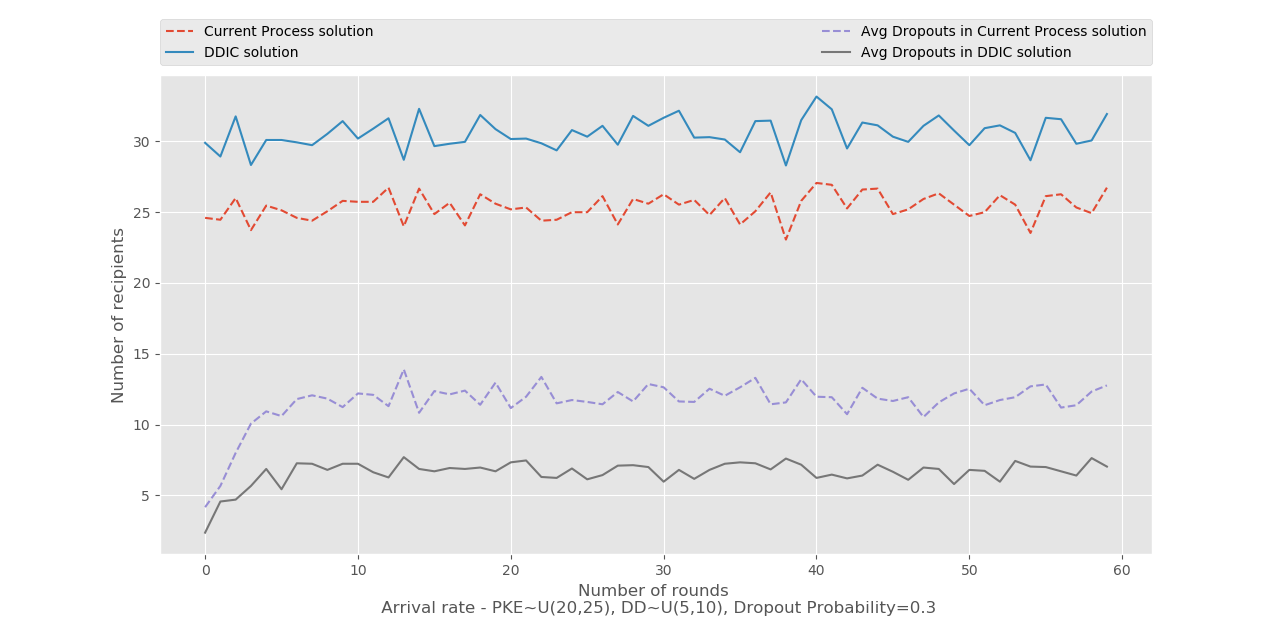}
	\end{minipage}
\end{figure}
\begin{figure}[H]
\hspace{-2.5cm}	\begin{minipage}[b]{0.6\textwidth}
		\includegraphics[scale=0.3]{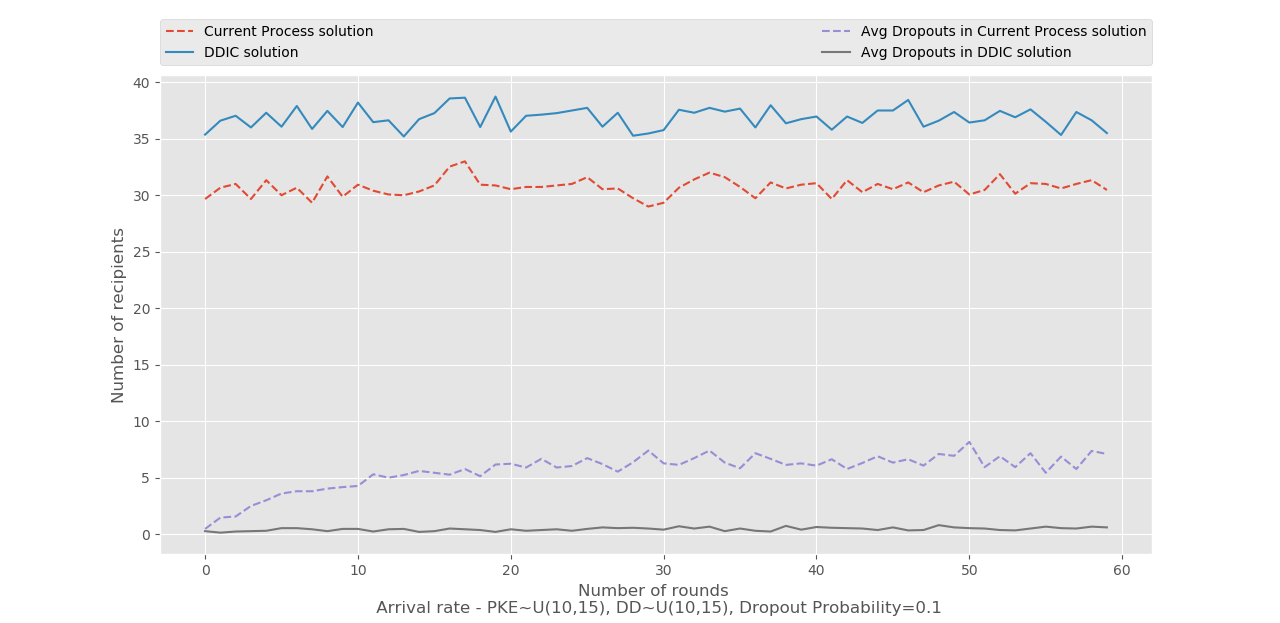}
	\end{minipage}
	\hspace{0.5cm}
	\begin{minipage}[b]{0.6\textwidth}
		\includegraphics[scale=0.3]{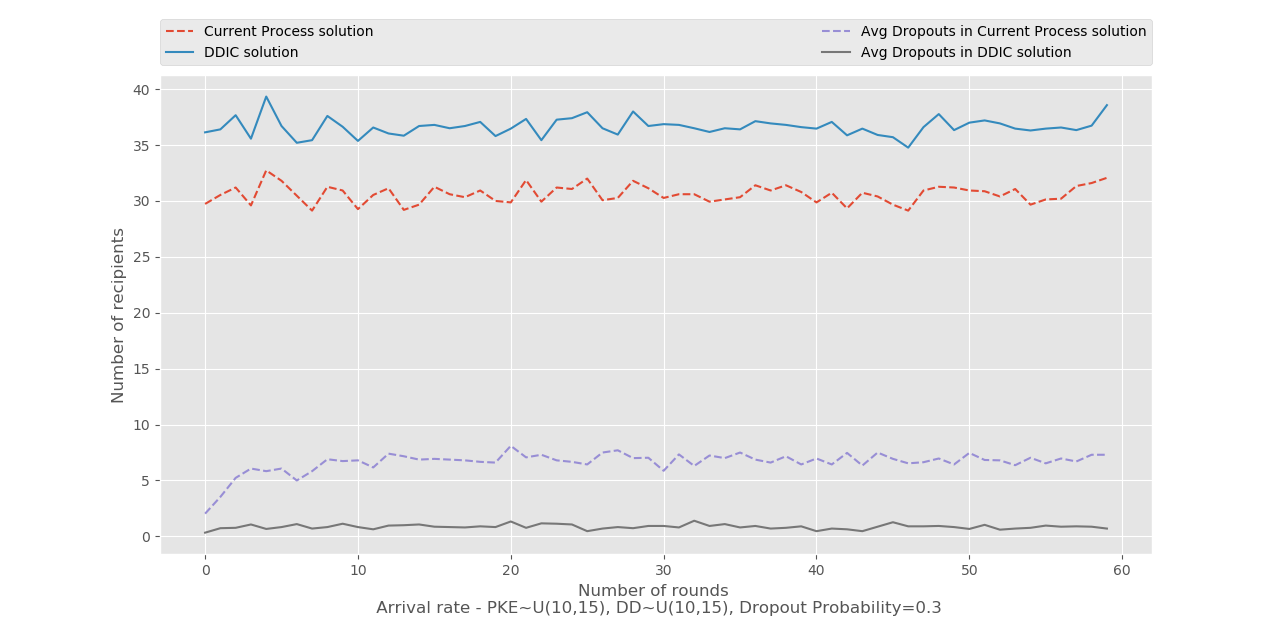}
	\end{minipage}
	
\hspace{-2.5cm}	\begin{minipage}[b]{0.6\textwidth}
		\includegraphics[scale=0.3]{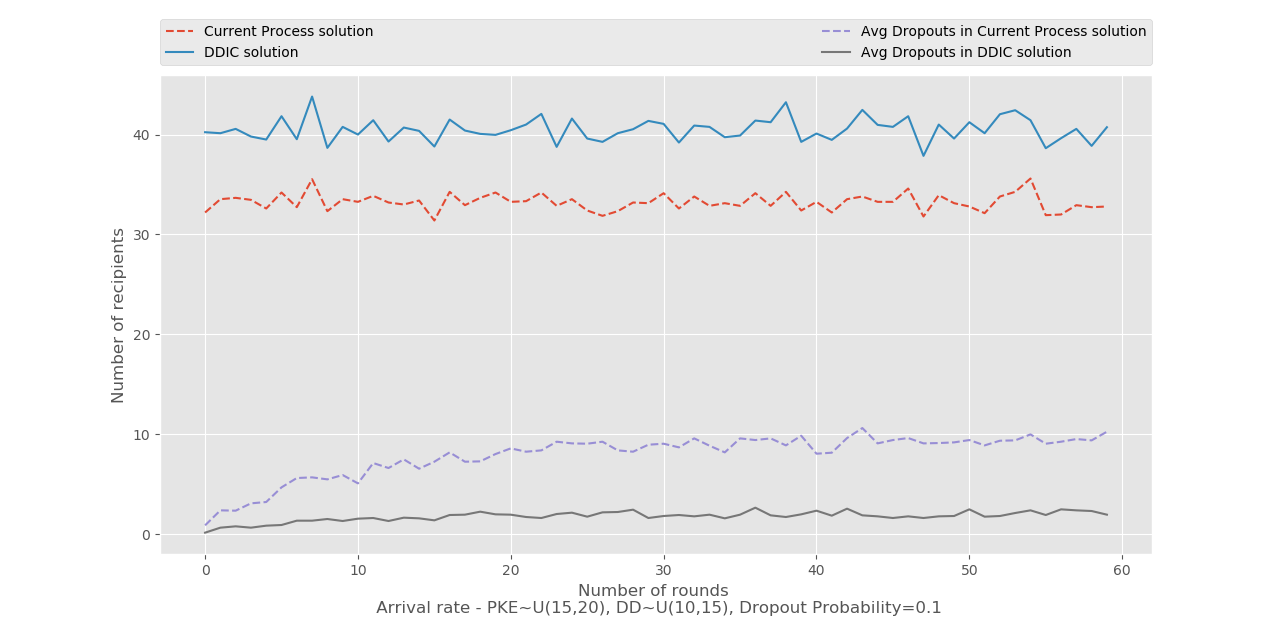}
	\end{minipage}
	\hspace{0.5cm}
	\begin{minipage}[b]{0.6\textwidth}
		\includegraphics[scale=0.3]{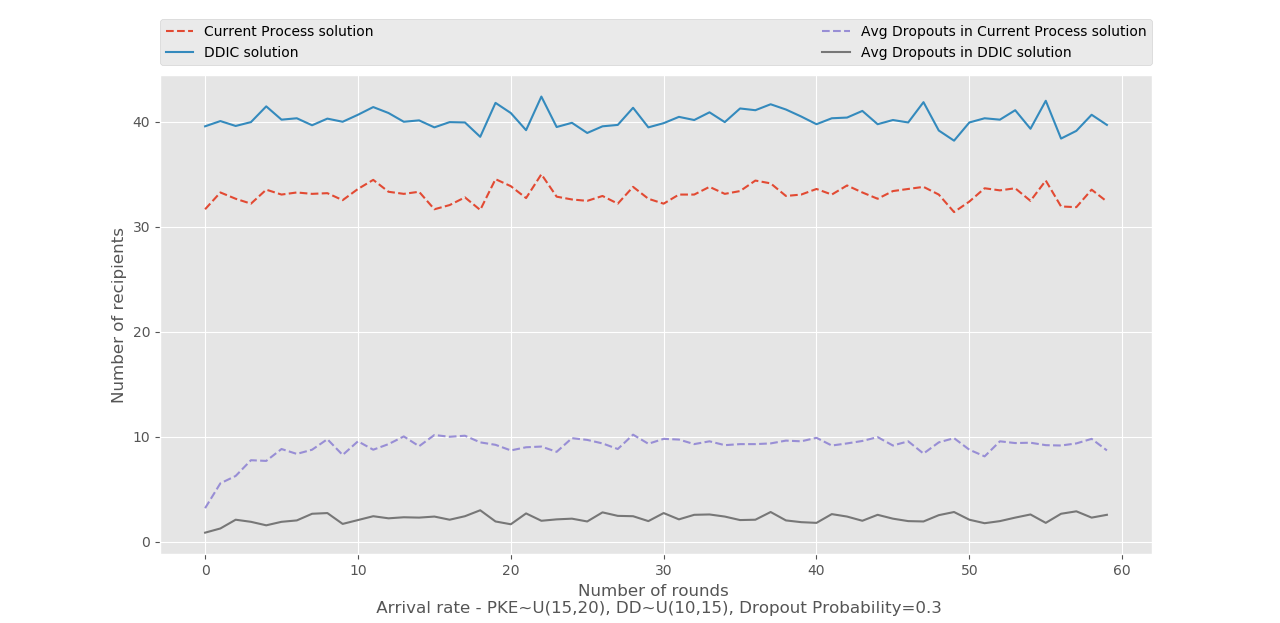}
	\end{minipage}
	
\hspace{-2.5cm}	\begin{minipage}[b]{0.6\textwidth}
		\includegraphics[scale=0.3]{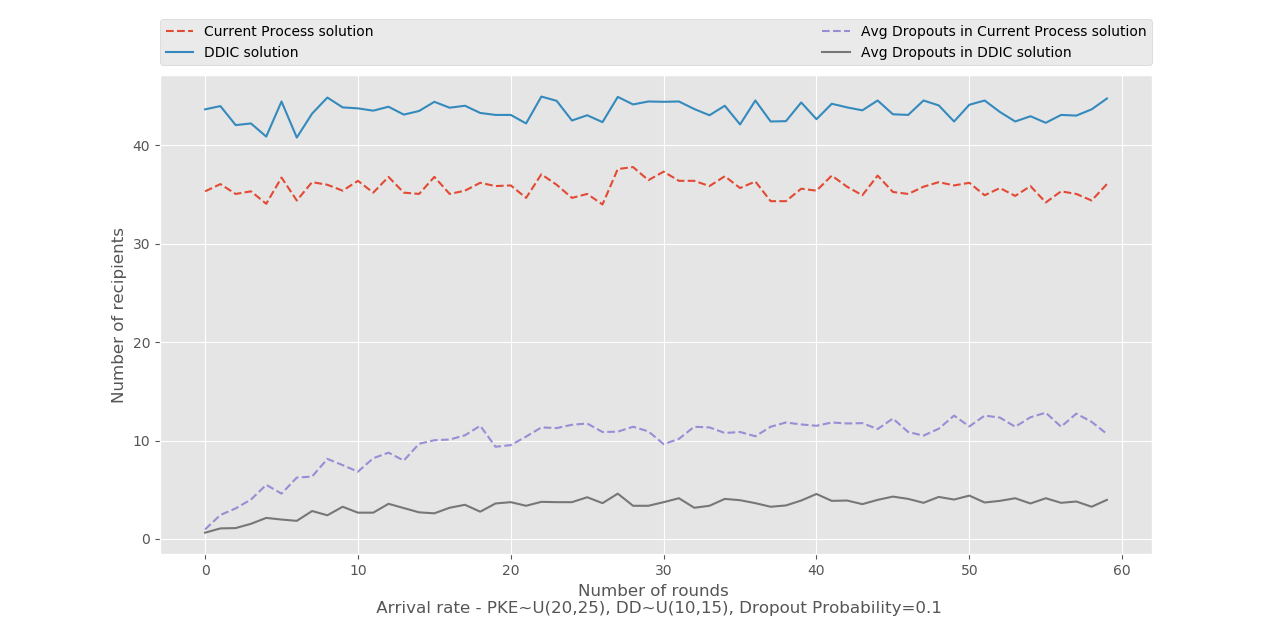}
	\end{minipage}
	\hspace{0.5cm}
	\begin{minipage}[b]{0.6\textwidth}
		\includegraphics[scale=0.3]{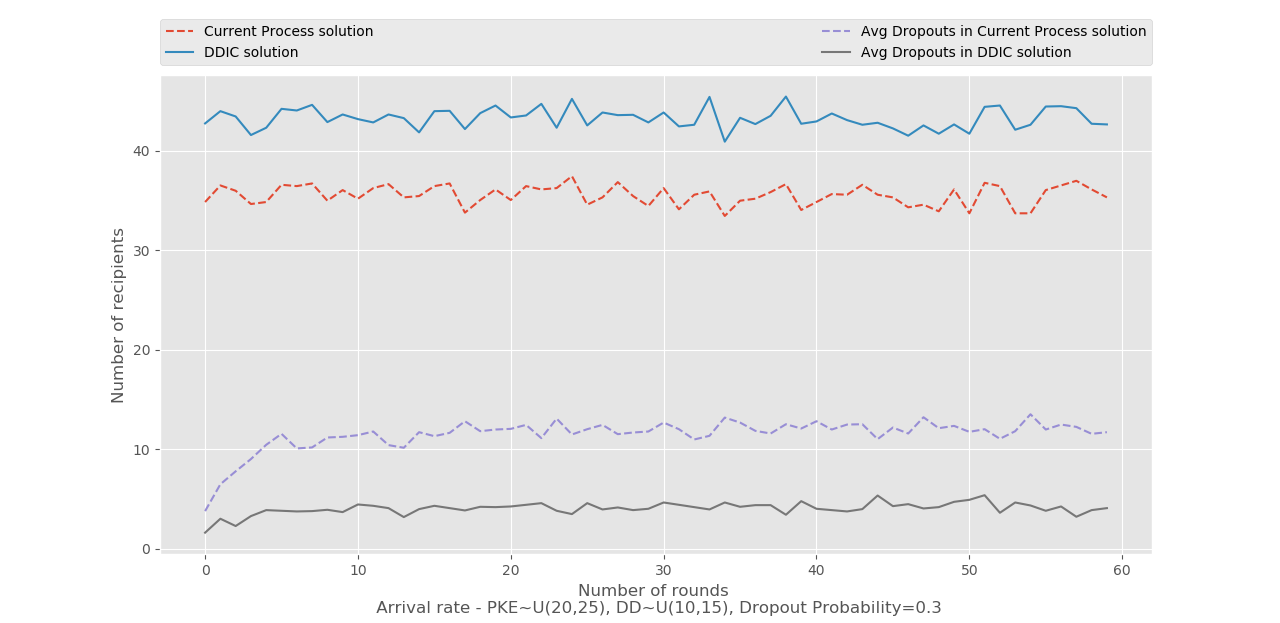}
	\end{minipage}
	
	\caption{Comparison of DDIC with current process in terms of number of patients matched and number of dropouts in each round, 18 instances have been shown with different arrival rates and dropout probabilities}
\end{figure}

\subsection{Blood group wise comparisons of waiting time and number of dropouts}

Since DDIC increases the possibility of compatible exchanges, waiting time becomes an important parameter to compare the output of both processes. We compared the average waiting time for both cases with different dropouts probabilities and results are shown in Table 2. 

\begin{table}
\caption{Average waiting time (in months) comparison between DDIC and Current Process (CP) for different blood groups in KEP registry. Note that: Waiting time for AB type patients in all cases were zero}

\resizebox{\textwidth}{!}{%
\begin{tabular}{p{5.5cm}c c c c c }
\hline
 & & O-type PT &  A-type PT & B-type PT \\
\hline
Arrival rates(/month) & DP & DDIC/CP & DDIC/CP & DDIC/CP \\
\hline
KEP=U(10-15),
DD=U(1-5) & 0.0 & 20.66/27.49 & 3.60/9.00 & 2.06/8.09 \\
 & 0.1& 5.31/7.04 & 1.08/2.27 & 0.68/1.85 \\
 & 0.3& 1.53/2.04 & 0.40/0.77 & 0.32/0.64 \\
\hline
KEP=U(15-20),
DD=U(1-5) & 0.0 & 22.60/27.44 & 4.90/8.78 & 2.51/7.21 \\
 & 0.1& 5.74/7.00 & 1.34/2.19 & 0.88/1.96 \\
 & 0.3& 1.68/2.03 & 0.47/0.70 & 0.31/0.56 \\
\hline
KEP=U(20-25),
DD=U(1-5) & 0.0 & 23.52/27.37 & 5.51/8.64 & 3.36/7.13 \\
 & 0.1& 6.03/6.95 & 1.57/2.25 & 0.94/1.93 \\
 & 0.3& 1.75/2.02 & 0.48/0.70 & 0.34/0.58 \\
\hline
KEP=U(10-15),
DD=U(5-10) & 0.0 & 11.48/27.86 & 0.52/9.00 & 0.27/7.47 \\
 & 0.1& 2.77/6.82 & 0.36/2.17 & 0.20/1.84 \\
 & 0.3& 0.85/2.06 & 0.20/0.73 & 0.11/0.61 \\
\hline
KEP=U(15-20),
DD=U(5-10) & 0.0 & 15.22/27.32 & 1.02/9.08 & 0.28/6.81 \\
 & 0.1& 3.98/6.95 & 0.44/2.17 & 0.27/1.99 \\
 & 0.3& 1.17/2.03 & 0.23/0.69 & 0.14/0.63 \\
\hline
KEP=U(20-25),
DD=U(5-10) & 0.0 & 18.27/27.49 & 1.93/8.93 & 0.82/7.49 \\
 & 0.1& 4.64/6.90 & 0.57/2.29 & 0.33/2.00 \\
 & 0.3& 1.34/2.01 & 0.26/0.66 & 0.16/0.59 \\
\hline
KEP=U(10-15),DD=U(10-15) & 0.0 & 1.97/27.49 & 0.16/8.85 & 0.06/7.47 \\
 & 0.1& 0.66/6.98 & 0.12/2.33 & 0.05/2.03 \\
 & 0.3& 0.31/2.03 & 0.07/0.74 & 0.04/0.62 \\
\hline
KEP=U(15-20),DD=U(10-15) & 0.0 & 7.53/27.52 & 0.18/8.69 & 0.12/8.26 \\
 & 0.1& 2.02/6.94 & 0.18/2.32 & 0.06/1.89 \\
 & 0.3& 0.62/2.04 & 0.10/0.71 & 0.05/0.60 \\
\hline
KEP=U(20-25),DD=U(10-15) & 0.0 & 11.59/27.33 & 0.32/8.89 & 0.13/8.01 \\
 & 0.1& 2.97/6.96 & 0.20/2.17 & 0.10/1.76 \\
 & 0.3& 0.88/2.00 & 0.12/0.68 & 0.07/0.60\\
\hline
\end{tabular}}
\end{table}

In Table 2, it can be observed that waiting time has reduced significantly for each blood group in the comparison. This reduction was expected as the number of recipients matched was higher in DDIC. Also, as the number of DD increases, the waiting time for each blood group decreases. This is due to the increase in the supply of kidneys to the KEP registry. Average waiting time for those who remain in KEP registry decreases with an increase in dropout probability, and that gives a more realistic scenario as pairs drop out from the KEP registry for various reasons.

\begin{table}
\caption{Comparison of total number of dropouts in DDIC and Current Process (CP) for different blood groups in KEP over 60 months registry(Replicated 30 times). Note that: Total number of dropouts for AB type patients in all cases were zero}

\resizebox{\textwidth}{!}{%
\begin{tabular}{p{5.5cm}c c c c c}
\hline
 & & O-type PT &  A-type PT & B-type PT \\
\hline
Arrival rates(/month) & DP & DDIC/CP & DDIC/CP & DDIC/CP \\
\hline
KEP=U(10-15),
DD=U(1-5)  & 0.1& 184.7/243.1 & 24.4/53.67 & 15.37/44.00 \\
 & 0.3& 208.17/272.23 & 39.3/66.63 & 29.9/58.3 \\
\hline
KEP=U(15-20),
DD=U(1-5) & 0.1& 281.03/337.5 & 41.87/74.2 & 29.03/63.11 \\
 & 0.3& 316.17/380.53 & 57.8/92.1 & 38.3/71.77 \\
\hline
KEP=U(20-25),
DD=U(1-5) & 0.1& 375.03/431.17 & 64.33/97.73 & 41.67/81.53 \\
 & 0.3& 423.63/487.10 & 78.73/115.30 & 58.3/93.6 \\
\hline
KEP=U(10-15),
DD=U(5-10) & 0.1& 98.37/240.37 & 8.06/54.5 & 4.8/44.87 \\
 & 0.3& 116.1/273.43 & 18.83/69.13 & 10.07/56.37 \\
\hline
KEP=U(15-20),
DD=U(5-10) & 0.1& 191.33/336.7 & 14.37/74.83 & 9.0/66.73 \\
 & 0.3& 217.87/378.4 & 27.6/87.8 & 19.0/81.93 \\
\hline
KEP=U(20-25),
DD=U(5-10) & 0.1& 287.17/429.19 & 24.33/95.7 & 14.1/82.3 \\
 & 0.3& 327.97/487.73 & 42.4/112.73 & 27.43/97.37 \\
\hline
KEP=U(10-15),DD=U(10-15) & 0.1& 22.8/238.53 & 2.9/54.9 & 1/47.4 \\
 & 0.3& 40.77/272.3 & 7.37/67.5 & 3.2/59.1 \\
\hline
KEP=U(15-20),DD=U(10-15) & 0.1& 99.23/339.87 & 5.83/76.43 & 2.2/60.33 \\
 & 0.3& 113.73/376.47 & 12.6/91.27 & 7.1/75.9 \\
\hline
KEP=U(20-25),DD=U(10-15) & 0.1& 187.57/428.4 & 9.2/93.13 & 5.13/74.6 \\
 & 0.3& 211.83/485.97 & 20.6/110.2 & 11.3/94.63 \\
\hline
\end{tabular}}
\end{table}

In table3, a comparison of the average number of dropouts over 60 rounds for each blood group is shown for both processes. It can be observed that each blood group benefits with DDIC in terms of the total number of dropouts as well. In scenarios with higher dropout probability, the relative benefit of DDIC increases.  

A retrospective study of deceased donor initiated chains is done by Cornelio et al. [9]. They have proposed a sequential procedure to create DD chains and worked on the historical data of Padua’s transplantation center. Here they have considered non-simultaneous executions of chains, so they try to find the longest possible chain through DD kidney. In India, non-simultaneous executions are not allowed, so that restricts the chain length. We did our simulation with the smallest possible chain (i.e. k=2). Our research aimed to observe the behavior of the DDIC mechanism over time and its effect on different parameters of the outcome. 

\section{Conclusion and Discussion}
Shortages of kidney donors have been a significant concern in providing compatible kidney transplants to needy recipients. KEPs were developed to increase the supply of compatible kidney to the pool of incompatible donor recipients pairs. With an increasing kidney exchange pool in India, not every recipient gets a compatible kidney in good time. The waiting time to get a compatible kidney in the KEP registry is increasing. In the last few years, deceased donation has increased steadily, and it is expected that with an increase in awareness about deceased donation, it will increase significantly in the country. With this prospect, the idea of DD initiated chains could increase the supply of compatible kidneys to the KEP registry, and also improve waiting time to get a compatible kidney, in both KEP as well as DD wait-list registry. An integer programming model is proposed in this paper, which allows us to create DD initiated chains. A long term simulation study is also done to analyze the gain that can be achieved through this. Results show that a significant gain in terms of the number of transplants and waiting time can be achieved for each blood group. There are challenges in the implementation of these DDIC simultaneously, as the DD kidney should be transplanted with 24 hours. In the Indian context, 5-way and 6-way simultaneous kidney exchange transplants have been successfully implemented [29,30], and as the transplant facilities increases in India, it is possible to conduct small simultaneous DDICs. The acceptance of a DD kidney in return for a living donor kidney is subjective, but as the demand for a kidney is very high, patients are willing to accept any possible transplant opportunity. 

In conclusion, DDIC seems to be one way forward to increase the utility of DD kidneys and the long term analysis shows that there will be a significant benefit in all the parameters in comparison to the current process.

\newpage

\section{Acknowledgements}

I would like thank Dr. Viswanath Billa and Dr. Deepa Usulumarty for representative data, discussion support and feedback.  

I would also like to thank Dr. Peter Biro for his suggestions and advice.

\end{document}